\documentclass[11pt,a4,leqno]{amsart}
\usepackage{latexsym,amsmath, amsfonts, graphics, epsf, epic}
\usepackage{amsthm}
\usepackage{amssymb}
\usepackage{amsopn}
\usepackage{amscd}
\usepackage{color}

\newcommand{\di}{\underline{i}}

\theoremstyle{plain}
\newtheorem{thm}{Theorem}[section]
\newtheorem*{thm*}{Theorem}  
\newtheorem*{idemthm*}{Idempotent Theorem}  
\newtheorem{lem}[thm]{Lemma}
\newtheorem*{lem*}{Lemma}    
\newtheorem*{orthlem*}{Orthogonality Lemma} 
\newtheorem{prop}[thm]{Proposition}
\newtheorem{cor}[thm]{Corollary}
\theoremstyle{definition}
\newtheorem{rk}[thm]{Remark}
\newtheorem*{rk*}{Remark}
\newtheorem{ex}[thm]{Example}
\newtheorem*{ex*}{Example}
\allowdisplaybreaks

\newcommand{\Z}{{\mathbb Z}}
\newcommand{\Q}{{\mathbb Q}}

\newcommand{\End}{\operatorname{End}}

\newcommand{\GL}{\mathsf{GL}}

\newcommand{\gl}{\mathfrak{gl}}
\renewcommand{\sl}{\mathfrak{sl}}

\newcommand{\divided}[2]{#1^{(#2)}}

\newcommand{\U}{\mathfrak{U}}

\title[Endomorphism rings of permutation modules]
{Endomorphism rings of permutation modules over maximal Young subgroups}
\author{Stephen Doty, Karin Erdmann and Anne Henke}
\date{\today}  
\thanks{All authors gratefully acknowledge support from Mathematisches 
Forschungsinstitut Oberwolfach, Research-in-Pairs Program. The second and
third author also thank the Bernoulli Center Lausanne.}

\begin{document}

\begin{abstract}
  Let $K$ be a field of characteristic two, and let $\lambda$ be a
  two-part partition of some natural number $r$. Denote the
  permutation module corresponding to the (maximal) Young subgroup
  $\Sigma_\lambda$ in $\Sigma_r$ by $M^\lambda$.  We construct a full
  set of orthogonal primitive idempotents of the centraliser
  subalgebra $S_K(\lambda) = 1_\lambda S_K(2,r)
  1_\lambda=\End_{K\Sigma_r}(M^\lambda)$ of the Schur algebra
  $S_K(2,r)$. These idempotents are naturally in one-to-one
  correspondence with the $2$-Kostka numbers.
\end{abstract}
\maketitle

\parindent=0pt


\section{Introduction}

Objects of central interest in the representation theory of symmetric
groups are permutation modules coming from actions on set
partitions. They provide a natural link with the representation theory
of general linear groups, via Schur algebras.
Assume $K$ is a field of positive characteristic $p$.  Fix natural
numbers $n$ and $r$ and fix partitions $\lambda$ and $\mu$ of $r$ of
not more than $n$-parts. The permutation module $M^\lambda$ over
$\Sigma_r$ is the module obtained by inducing the trivial
representation from the Young subgroup $\Sigma_\lambda$ to the
symmetric group $\Sigma_r$.  The indecomposable direct summands of
$M^{\lambda}$ are known as Young modules.  By James' submodule theorem
\cite[7.1.7]{JK} there is a unique indecomposable summand of
$M^{\lambda}$ containing the Specht module $S^{\lambda}$. This summand
is by definition the Young module $Y^{\lambda}$.  The module
$M^{\lambda}$ is in general a direct sum of Young modules $Y^\mu$, and
if $Y^{\mu}$ occurs as a summand then $\mu \geq \lambda$, in the usual
dominance order on partitions.  The $p$-Kostka number $[M^\lambda :
Y^\mu]$ is the number of indecomposable summands of $M^\lambda$
isomorphic with $Y^\mu$.  Thus we have:
\begin{equation*} \textstyle
M^\lambda \simeq  
\bigoplus_{\mu \geq \lambda} \;\; [M^\lambda : Y^\mu]\;  Y^\mu .
\end{equation*}
Note that if $\lambda$ is an $n$-part  composition 
of $r$ then we still have the permutation module
$M^\lambda$ defined as above. In this situation  
if $\lambda_0$ is the partition obtained from $\lambda$ by 
ordering its parts then 
$M^{\lambda_0}
\simeq M^\lambda$.

Let $S_K(n,r)$ be the Schur algebra of degree $r$, 
then  $S_K(n,r)$ is given by
\begin{equation*}
S_K(n,r) = \End_{K\Sigma_r}(E^{\otimes r}) \simeq \textstyle
\End_{K\Sigma_r}(\bigoplus_{\lambda} M^\lambda) 
\end{equation*}
where $E$ is a given $n$-dimensional $K$-vector space, and $\lambda$ 
varies  over the set $\Lambda(n,r)$ of $n$-part compositions of $r$.
For the connection of Schur algebras with general linear groups, see
Green~\cite{green}.  The idempotent $1_\lambda \in S_K(n,r)$ corresponds
to the projection onto $M^{\lambda}$ with kernel $\oplus _{\mu \neq
\lambda} M^{\mu}$. We define the centraliser subalgebra $S_K(\lambda)$
of $S_K(n,r)$ by
\begin{equation*}
S_K(\lambda) = 1_\lambda S_K(n,r) 1_\lambda \simeq 
\End_{K\Sigma_r}(M^\lambda).
\end{equation*}
In this paper we study these algebras when $\lambda$ is a partition of
at most two parts; that is, the associated Young subgroup is maximal.
Then it is known (see for example \cite{james}, Example 14.4) that the
ordinary character of $M^\lambda$ is multiplicity-free.  It follows
that the algebra $S_K(\lambda)$ is commutative (see \cite{scott}, and
see also Remark \ref{KE:rmk} below), and that any given Young module
$Y^{\mu}$ occurs at most once as a direct summand of $M^{\lambda}$ .
All idempotents of $S_K(\lambda)$ are central, and there are finitely
many primitive idempotents, in one-to-one correspondence with the
indecomposable summands of $M^{\lambda}$. The blocks of $S_K(\lambda)$
are therefore precisely the endomorphism rings of the Young modules
$Y^{\mu}$ which occur as a direct summand of $M^{\lambda}$.

Our main result is an explicit construction of a full set of
orthogonal primitive idempotents of the algebra $S_K(\lambda)$ where
$\lambda$ is a partition of at most two parts and char$(K)=2$. These
idempotents are naturally in one-to-one correspondence with the
$2$-Kostka numbers.  The philosophy is to consider an infinite family
of algebras at the same time, as was done in~\cite{DEH1}. This is
possible, by exploiting the presentation obtained in \cite{DG:PSA} of
the Schur algebra as quotient of the universal enveloping algebra.
This approach allows one to keep $m=\lambda_1-\lambda_2$ fixed and let
$r = \lambda_1+\lambda_2$ vary arbitrarily.
Our motivation is to describe idempotents explicitly; this
is a notoriously hard problem, in general, especially in 
the modular setting. We solve this problem completely, for
our situation, when $p=2$; the case of odd primes seems to be more complicated.
The results on idempotents
in this paper can be thought of as an algebraic realization  of the
combinatorial description of the quarter-infinite Kostka matrix in Section 2.1
below. 
%

\section{Main Results} \label{sec2}

\subsection{The $p$-Kostka matrix} 
We fix some notation.  Let $K$ be a field of positive characteristic
$p$. For any natural number $r$, we let $\lambda=(r-k,k)$ and
$\mu=(r-s,s)$ vary over the two-part partitions of $r$. The $p$-Kostka
numbers $[M^{(r-k,k)}:Y^{(r-s,s)}]$ do not depend on $r$ but only on
$m:= \lambda_1 - \lambda_2=r-2k$ and $g:= \lambda_2- \mu_2= k-s$.  So
they can be described by a quarter-infinite matrix with
$(m,g)$th-entry the above p-Kostka number.  Set
\begin{equation*}
B(m,g)=\tbinom{m+2g}{g}.
\end{equation*}
By \cite{He1, henke} it is known that $Y^{(r-s,s)}$ is a direct
summand of $M^{(r-k,k)}$ if and only if $B(m,g)=\binom{r-2s}{k-s}\neq
0$ modulo $p$.  Since the multiplicity $[M^{(r-k,k)}: Y^{(r-s,s)}]$ is
at most one, the $(m,g)$th entry of the $p$-Kostka matrix is one if
$B(m,g) \ne 0$ modulo $p$ and zero otherwise. This latter result is
based on a general formula by Klyachko~\cite{klyachko}, Corollary 9.2,
reformulated by Donkin~\cite{donkin} in (3.6). However, neither 
reference gives an explicit answer.

\subsection{Notation}
We need the $p$-adic expansion of integers. If $a=\sum^s_{j=0} a_jp^j$
with $0\leq a_j\leq p-1$ for all $j$ then we will write $a = [a_0, a_1,
\ldots, a_s]$.
It is a well-known property of binomial coefficients that 
\begin{equation*}
B(m,g) \equiv \prod_i \tbinom{(m+2g)_i}{g_i} \mod{p}.
\end{equation*}
We call $\prod_i \binom{(m+2g)_i}{g_i}$ modulo $p$ the binomial expansion
of $B(m,g)$ and we write $B(m,g)_i$ for the $i$-th factor of this
product. Sometimes we will also write the binomial coefficient modulo
$p$ as a matrix with two rows where the $i$-th column is the $i$-th
factor of the product, for $i \geq 0$:
\begin{equation*}
B(m,g)= \left( 
\begin{smallmatrix}
(m+2g)_0 &(m+2g)_1 & \ldots &(m+2g)_i & \ldots\\
g_0 &g_1 & \ldots & g_i & \ldots 
\end{smallmatrix}\right).
\end{equation*}

\subsection{The canonical basis of $S_K(\lambda)$}
Let $m\geq 0$. We consider the infinite family of algebras
$S_K(\lambda)$ where $\lambda$ runs through all partitions $\lambda =
(r-k, k)$ such that $r-2k=m$. The presentation from \cite{DG:PSA}
provides $S_K(\lambda)$ with a basis $\{ b(a): 0\leq a\leq k\}$ with
good properties.  In fact, it is not difficult to see that this basis
is inherited from the canonical basis of Lusztig's modified form
$\dot{\U}(\sl_2)$ of the enveloping algebra $\U(\sl_2)$, which is
worked out as an example in \cite[\S29.4.3]{Lusztig}.  In particular,
the product $b(i)b(j)$ depends only on $m$ (and not on the degree
$r$), where terms $b(s)$ appearing with $s>k$ are set zero. See \S3
below for more details on this basis, upon which our computation of
primitive idempotents is based.

One can also consider the infinite-dimensional generic algebra
$\dot{\U}(\lambda) = 1_\lambda \dot{\U}(\gl_2) 1_\lambda$, as in
\cite{DEH1}, which embeds naturally in the inverse limit of the
following sequence of surjections 
\begin{equation*}
  \cdots S_\Q(\lambda+\delta) \to S_\Q(\lambda) \to \cdots \to
  S_\Q(\nu+\delta) \to S_\Q(\nu)
\end{equation*}
where $\delta = (1,1)$ and $\nu = (r-2k,0)=(m,0)$. Here the maps are
induced from corresponding maps on the Schur algebra level,
corresponding to tensoring by the determinant representation.  (More
precisely, one tensors a corresponding coalgebra by the determinant to
get an embedding of coalgebras, and then dualizes to get s surjection
between Schur algebras.)

The algebra $S_\Q(\lambda)$ is a homomorphic image of
$\dot{\U}(\lambda)$ and our multiplication formula for the canonical
basis elements of $S_K(\lambda)$ is coming from a corresponding
multiplication formula in the generic algebra. (The setup is
compatible with change of base ring.)  The generic point of view is
closely related to the approach of \cite{BLM}; see \cite{DEH1} for
further details.

\subsection{Construction of primitive idempotents} \label{sub:main}
We now take char$(K)=2$ and keep $m\geq 0$ fixed. We work in an
algebra $S_K(\lambda)$ of large enough degree $r$ (of the right
parity).
For any $m\geq 0$ and $g \geq 0$ such that $B(m,g)$ is non-zero modulo
two and the degree $r$ is large enough (that is ~$r \geq m+2g$) we
will now define elements in the algebra $S_K(\lambda)$. First we
introduce two index sets: let
\begin{equation*}
\begin{aligned}
I_{m,g}&:= \{ u: g_u=0 \mbox{ and } (m+2g)_u=1 \}, \\ 
J_{m,g}&:= \{ u: g_u=1 \mbox{ and } (m+2g)_u=1 \}.
\end{aligned}
\end{equation*}
Then for a natural number $t$ define elements in the algebra
$S_K(\lambda)$ by
\begin{eqnarray}\label{eqn100}
  e_{m,g}&:= &\prod_{u\in J_{m,g}} b(2^u) \prod_{u\in I_{m,g}}
  (1-b(2^u)). \\
  (e_{m,g})_{\leq t}&:= &\prod_{u\in J_{m,g}, u \leq t} b(2^u) \prod_{u\in I_{m,g}, u\leq t}
  (1-b(2^u)). \nonumber
\end{eqnarray}
\begin{rk*}
We can associate to each factor of the binary expansion
of $B(m,g)$ a factor of an element $e_{m,g}$ by the following rule:
\begin{center} \renewcommand{\arraystretch}{1.5} 
\begin{tabular}{|c|c|c|c|c|}\hline
$B(m,g)_u$ & $\binom{1}{1}$ & $\binom{1}{0}$ & $\binom{0}{0}$ 
           & $\binom{0}{1}$\\ \hline 
Factor of $e_{m,g}$ & $b(2^u)$ & $(1-b(2^u))$ & $1$ & $0$ \\ \hline
\end{tabular}
\end{center}
In particular, an element $e_{m,g}$ defined in this way would be zero
if and only if $\tbinom{0}{1}$ occurs in the binary expansion of
$B(m,g)$, that is if and only if $B(m,g)= 0$ modulo two.
\end{rk*}

The main result of the paper is the following:

\begin{idemthm*} 
For any fixed $m \geq 0$, the set of elements $e_{m, g}$ with $B(m,g)
\neq 0$ modulo two and $m+2g \leq r$ is a complete set of primitive
orthogonal idempotents for the algebra $S_K(\lambda)$.
\end{idemthm*}


This theorem will be proved at the end of Section~\ref{sec6}. In fact
parts of the proof of this result are not so difficult to see.
Observe the following:
\begin{itemize}
\item[(i)] { The element $e_{m,g}$ is non-zero.}  By
  Proposition~\ref{prop2} or Lemma~\ref{lem12} one can even express it
  explicitly as a linear combination of the basis elements.
\item[(ii)] { If $g\neq d$ and $B(m,g)$ and $B(m,d)$ are both
    non-zero modulo two then $$e_{m,g}^2\cdot e_{m,d}^2=0.$$}
\begin{proof}[Proof of (ii)]
Let $i$ be minimal such that $B(m,g)_i\neq B(m,d)_i$.  Since columns
$<i$ are the same, and both binomial coefficients are non-zero, the
i-th columns cannot be zero: Suppose that one is zero, the other
not. Then $(m+2d)_i \not \equiv (m+2g)_i$. However, in column $i$ the
carry overs from the previous columns are the same, say $x$, and
$d_{i-1}=g_{i-1}$. This implies a contradiction:
\[
(m+2d)_i = m_i+ d_{i-1} + x =m_i+g_{i-1}+x = (m+2g)_i  \mod{2}.
\]
Hence one of them is $\tbinom{1}{1}$ and the other is $\tbinom{1}{0}$.
So the squares of the elements in the algebra have factors $b(2^i)^2$
and $(1-b(2^i)^2)$ respectively.  In Section~\ref{sec4} we will show
that the elements $b(2^i)^2$ are idempotents (see
Proposition~\ref{cor2}), and this implies that the product is zero.
\end{proof}
\end{itemize}

Furthermore the number of primitive idempotents of $S_K(\lambda)$ is
equal to the number of non-zero binomial coefficients,
by~\cite{henke}.  So, once we have established that the elements in
question are idempotents then the theorem is proved. The fact that the
elements in question are idempotents will follow from an orthogonality
result:

\begin{orthlem*} 
Suppose $B(m,g)_s$ is zero, then $e_{m,g}^2\cdot b(2^s)^2=0$.
\end{orthlem*}

The proof of the Orthogonality Lemma is given in Section~\ref{orthlem}.

\subsection{Blocks of the algebra $S_K(\lambda)$}
Recall that a block of a finite-dimensional algebra $A$ is given by
$eA$ where $e$ is a central idempotent which is primitive, viewed as
an element of the centre of $A$.  By the Idempotent Theorem, a block
of the algebra $S_K(\lambda)$ has the form $e_{m,g} S_K(\lambda)$.
The block has basis
\begin{equation*}
\{ e_{m,g}b(a)\ :  a = [a_0, a_1, \ldots ] \mbox{ where $a_s=1$ only for }
(m+2g)_s=0\}.
\end{equation*}
To see that this set is linearly independent one uses Lemma
\ref{lem12}; and to show that it spans, one observs using Lemma
\ref{lem13} and Lemma \ref{lem2} that if $(m+2g)_s\neq 0$ then
$e_{m,g}b(2^s)$ can be expressed in terms of the given set, by
elements of 'lower degree'.  By Lemma~\ref{lem12}, the block is
(minimally) generated as an algebra by all
\begin{equation*}
\{ e_{m,g}b(2^s): \ s\geq 0 \mbox{ and } (m+2g)_s=0\}.
\end{equation*}
By the Orthogonality Lemma, the block has a set of generators
with square zero. Hence for a general degree $r$, this block is
isomorphic to a quotient of an algebra of the form
\begin{equation*} \textstyle
\bigotimes K[x_i]/\langle x_i^2\rangle.
\end{equation*}
a tensor product of finitely many local two dimensional algebras.


\section{Basis and  multiplication structure in $S_K(\lambda)$} 
\label{sec3}
In this section we take $\lambda=(\lambda_1, \lambda_2)$ to be a
two-part partition and we study the multiplicative structure of
$S_K(\lambda)$ over a field $K$ of characteristic $p \geq 0$.  The
results from this section will then be used to obtain in
characteristic two a reduction formula for $b(2^s)^2$ (see
Section~\ref{sec4}).

We describe briefly some results from \cite{DG:PSA0}, \cite{DG:PSA}.
Over $\Q$, the Schur algebra $S_\Q(2,r)$ is isomorphic to the quotient
of the universal enveloping algebra $\U(\gl_2)$ modulo the ideal
generated by
\begin{equation*}
H_1(H_1-1) \ldots \cdot (H_1-r);
\end{equation*}
alternatively, $S_\Q(2,r)$ can be described as the quotient of
$\U(\sl_2)$ modulo the ideal generated by
\begin{equation*}
(h+r)(h+r-2)\cdots (h-r+2)(h-r).
\end{equation*}
Here, as a basis for the Lie algebra ${\gl}_2$ one takes $e=e_{12},
f=e_{21}$ as usual and $H_1, H_2$ respectively the diagonal matrices
$e_{11}$ and $e_{22}$, where $e_{ij}$ is the usual matrix unit, and as
basis for $\sl_2$ the usual $e,f,h$ where $h=H_1-H_2$ is the
commutator of $e$ and $f$.

The family of algebras $\{ S_K(2,r) \}_K$ ($K$ a field) is defined
over $\Z$ using the usual divided powers.  In this presentation, the
idempotent $1_{\lambda}$ which we earlier defined as projection
corresponding to $\lambda$, is equal to the image (in the Schur
algebra) of
\begin{equation*}
1_{\lambda} =  
\tbinom{H_1}{\lambda_1}\tbinom{H_2}{\lambda_2};
\end{equation*} 
see \cite[Lemma 5.3]{DEH1}.  
Now $S_\Q(\lambda) = 1_\lambda S_\Q(2,r) 1_\lambda$. There is a
natural $\Z$-form $S_\Z(2,r)$ of $S_\Q(2,r)$, namely the image of the Kostant
$\Z$-form of $\U(\gl_2)$ (or of $\U(\sl_2)$) under the quotient map
$\U(\gl_2) \to S_\Q(2,r)$.  Set $S_\Z(\lambda) = S_\Q(\lambda) \cap
S_\Z(2,r)$.  If we set
\begin{equation*}
b(i):= 1_{\lambda}f^{(i)}e^{(i)}1_{\lambda} \in S_{\Q}(2,r)
\end{equation*} 
then $S_{\Z}(\lambda)$ is the subalgebra of $S_\Q(\lambda)$ with
basis $\{ b(0), b(1), \ldots, b(\lambda_2)\}$.  Here 
\begin{equation*} \label{eqn15}  
  e^{(i)}:= \tfrac{e^i}{i!} \quad \mbox{
    and } \quad f^{(i)}:= \tfrac{f^i}{i!}
\end{equation*} 
are the usual divided powers in the enveloping algebra.

We next describe the multiplicative structure of $S_\Z(\lambda)$.
Most important for us is a multiplication formula for the basis
elements $b(i)$, given in Proposition~\ref{prop2}, which also proves
again that the $b(i)$ generate a $\Z$-form of $S_\Q(\lambda)$. In what
follows, we identify generators $e,f,h, H_1, H_2$ with their images in
the quotient $S_\Q(\lambda)$.

The following formula, valid in $\U(\sl_2)$, is easily derived by
induction on $a$: 
\begin{equation} \label{eqn12}
e {f}^{a} = {f}^{a} e + a {f}^{a-1}(h-a+1).
\end{equation}
Since this formula holds in the enveloping algebra (over $\Q$), it is
valid in its homomorphic image $S_{\Q}(2,r)$.  The first part of the
next lemma is contained in \cite{DG:PSA}.

\begin{lem}  \label{lem14}
In $S_{\Q}(2,r)$ we have the equality $h 1_\lambda = m 1_\lambda$, where $m
= \lambda_1 -\lambda_2$. Moreover, for any $k$ we have
\[
b(1) \cdot b(k) = (k+1)^2 b({k+1})+ k(m+k+1) b(k).
\]
\end{lem}

\begin{proof}  
To see this, first calculate using formula (\ref{eqn12}):
\begin{align*}
(k!)^2 \cdot b(1) \cdot b(k) &= fef^ke^k 1_\lambda \\
 &= f(f^ke+kf^{k-1}(h-k+1))e^k 1_\lambda \\
 &= f^{k+1}e^{k+1}1_\lambda + kf^k(h-k+1)e^k 1_\lambda \\
 &= f^{k+1}e^{k+1}1_\lambda + kf^ke^k(h+k+1) 1_\lambda 
\end{align*}
where we have used the fact that $he^k = e^k(h+2k)$. This holds in the
enveloping algebra of $\mathfrak{gl}_2$, and hence is valid in
$S_{\Q}(2,r)$.  Now apply the first statement of this Lemma to obtain the
desired formula.
\end{proof}

\begin{lem} \label{lem11}
Set $x=b(1)$. Then we have in $S_{\Q}(2,r)$ for any $k \ge 1$ the
equality
\[
b({k+1})= \tfrac{1}{(k+1)!^2}x(x-(m+2))(x-2(m+3))\cdots (x-k(m+k+1)).
\]
\end{lem}

\begin{proof} Proceed by induction on $k$.  Define
\begin{equation*}
F_{k+1}(x)= x(x-(m+2))(x-2(m+3))\cdots (x-k(m+k+1)).
\end{equation*}
The case $k=1$
in the preceding lemma gives the equality
\begin{equation*}
b(2) = \tfrac{1}{2^2}(x^2 - (m+2)x)= \tfrac{F_2(x)}{(2!)^2}.
\end{equation*}
Thus the formula of the lemma is valid in case $k=1$.  Assume that
$b(k)=\frac{F_k(x)}{(k!)^2}$.
By the preceding lemma and the inductive hypothesis we then have
\begin{equation*}
\begin{aligned}
b({k+1}) &= \tfrac{1}{(k+1)^2} \cdot (b(1) b(k) - k(m+k+1) b(k)) \\ 
&=
\tfrac{1}{{(k+1)!}^2} \cdot (x-k(m+k+1)) F_{k}(x) =
\tfrac{1}{{(k+1)!}^2} F_{k+1}(x).
\end{aligned}
\end{equation*}
\end{proof}

\begin{prop} \label{prop3}
The algebra $S_\Q(\lambda)$ is semisimple and generated by $b(1)$.
\end{prop}

\begin{proof}
  The semisimplicity statement is clear, since $M^\lambda$ is
  completely reducible as a $\Sigma_r$-module in characteristic zero.
  Thus only the claim about generation needs to be proved.  It follows
  from Lemma \ref{lem11} that we have the equality
\begin{equation*}
b(k) = 1_\lambda \divided{f}{k} \divided{e}{k} 1_\lambda 
= \tfrac{F_k(x)}{(k!)^2}
\end{equation*}
for all $k\ge 2$.  This formula holds in $S_\Q(2,r)$ and hence any
element in $S_\Q(\lambda)$ is generated by $x=b(1)$. 
\end{proof}

\begin{prop} \label{commprop}
  The algebra  $S_\Q(\lambda)$ is isomorphic with  $\Q[T]/
  (F_{\lambda_2+1}(T))$.
\end{prop}

\begin{proof} 
By commutation formulas appearing in \cite{DG:PSA} we have
\begin{equation*}
b({\lambda_2+1})=1_\lambda \divided{f}{\lambda_2+1}
\divided{e}{\lambda_2+1} 1_\lambda = 0 =F_{\lambda_2+1}(x)
\end{equation*}
since $\lambda + (\lambda_2+1)(1,-1) = (\lambda_1+\lambda_2+1, -1)$ is
not a polynomial weight belonging to $\Lambda(2,r)$, for any $\lambda$.
The proposition now follows from Lemma \ref{lem11}.
\end{proof}

\begin{rk}\label{KE:rmk} 
It follows immediately from the preceding Proposition 
that the algebra $S_\Q(\lambda)$ is a commutative algebra.
In fact, the commutativity of $S_\Q(\lambda)$ is a consequence of the
fact that the permutation module $M^\lambda$ is multiplicity-free (see
\cite{scott}). This $\Sigma_n$-module is semisimple and its
composition factors are absolutely irreducible, so by Schur's Lemma
$\End_{\Q\Sigma_n}(M^\lambda)$ is a direct sum of copies of the field
$\Q$.
\end{rk}

\begin{prop} \label{prop2} A multiplication formula for the basis
  elements is given by:
\begin{equation} \label{eqn1a}
b(i)\cdot b(j) = \sum_{k=0}^i \tbinom{j+k}{i} \tbinom{j+k}{k} \tbinom{m+j+i}
{i-k}\, b(j+k).
\end{equation}
\end{prop}      

When $a > \lambda_2$ then $b(a)$ is zero in this formula.

\begin{proof} The proof is by induction on $i$. The induction
  beginning for $i=1$ is given by Lemma~\ref{lem14}. Let now $i>1$.
  Then using Lemma~\ref{lem14}, the product $P:= b({i+1}) \cdot b(j)$
  equals:
\begin{eqnarray*}
P &=& \tfrac{b(i) \cdot (b(1)-i(m+i+1))}{(i+1)^2} \; b(j)  \\
&=&  \tfrac{b(1)-(j+k)(m+j+k+1) 
    + (j+k) (m+j+k+1) - i(m+i+1)}{(i+1)^2} \; b(i)\, b(j)\\ 
&=&  \tfrac{b(1)-(j+k)(m+j+k+1) 
    + (j+k-i) (m+j+k+i+1)}{(i+1)^2}  \; b(i)\, b(j) \\
&=& \sum_{k=0}^{i} \tfrac{b(1)-(j+k)(m+j+k+1))}{(i+1)^2} \tbinom{j+k}{i}
    \tbinom{j+k}{k} \tbinom{m+i+j}{i-k} \; b({j+k}) \\
&&  + \sum_{k=0}^{i} \tfrac{(j+k-i)(m+j+k+i+1))}{(i+1)^2} \tbinom{j+k}{i}
    \tbinom{j+k}{k} \tbinom{m+i+j}{i-k} \; b({j+k}) \\
&=& \sum_{k=0}^{i} \tfrac{k+1}{i+1} \tbinom{j+k+1}{i+1}
    \tbinom{j+k+1}{k+1} \tbinom{m+i+j}{i-k} \; b({j+k+1}) \\
&&  + \sum_{k=0}^{i} \tfrac{(m+j+k+i+1))}{(i+1)} \tbinom{j+k}{i+1} 
    \tbinom{j+k}{k} \tbinom{m+i+j}{i-k} \; b({j+k}) \\
&=& \sum_{k=1}^{i+1} \tfrac{k}{i+1} \tbinom{j+k}{i+1}
    \tbinom{j+k}{k} \tbinom{m+i+j}{i+1-k} \; b({j+k}) \\
&&  + \sum_{k=0}^{i} \tfrac{(m+j+k+i+1))}{(i+1)} \tbinom{j+k}{i+1}
    \tbinom{j+k}{k} \tbinom{m+i+j}{i-k} \; b({j+k}) \\
&=& \sum_{k=0}^{i+1} \tbinom{j+k}{i+1} \tbinom{j+k}{k} 
      \tbinom{m+j+i+1}{{i+1-k}} \cdot  b({j+k})
\end{eqnarray*}
and the induction is complete.
Note that the last equality above is justified as follows:
\begin{eqnarray*}
   \tfrac{k}{i+1} \tbinom{m+i+j}{i+1-k} 
      + \tfrac{m+j+k+i+1}{i+1} \tbinom{m+i+j}{i-k} 
  &=& \tfrac{k}{i+1} \tbinom{m+i+j}{i+1-k} + \tfrac{i+1-k}{i+1} 
      \tbinom{m+i+j}{i+1-k} + \tbinom{m+i+j}{i-k} \\
  &=& \tbinom{m+i+j}{i+1-k} + \tbinom{m+i+j}{i-k} \\
  &=& \tbinom{m+i+1+j}{i+1-k} 
\end{eqnarray*}
for any $k$ satisfying $1 \leq k \leq i$.
\end{proof}

The previous proposition verifies again that $S_\Z(\lambda)$ is a
$\Z$-form. We have $S_K(\lambda) \simeq S_{\Z}(\lambda) \otimes_\Z K$
for any field $K$, and by abuse of notation we still write $b(i)$ for
$1_\lambda \divided{f}{i} \divided{e}{i} 1_\lambda \otimes 1_K$.  Then
the multiplication formula (\ref{eqn1a}) is also valid in the
$K$-algebra $S_K(\lambda)$; and again $b(a)=0$ whenever $a> \lambda_2$.

For the remainder of this section, we work over the field $K$ of
positive characteristic $p$.


\begin{lem} \label{lem12}
  Write $i=[i_0, i_1, \ldots]$ $p$-adically. Then $b(i) = \prod_{t
  \geq 0} b({i_t \cdot p^{t}})$ in $S_K(\lambda)$.
\end{lem}

\begin{proof}
  This is shown by induction on the length of the $p$-adic
  decomposition of $i$. Assume that $j=[i_0, i_1, \ldots, i_{t-1}]$;
  then by Equation (\ref{eqn1a}):%
\begin{equation*}
b(j) \cdot b({i_t p^t}) = \sum_{k=0}^j  \tbinom{i_tp^t+k}{j}
\tbinom{i_tp^t+k}{k} \tbinom{m+i_tp^t+j}{j-k} \; b({i_t p^t +k}).
\end{equation*} 
In the above sum, the binomial coefficient $\tbinom{i_t p^t +k}{j}=
\tbinom{k}{j}$ is nonzero if and only if $j$ is $p$-contained in $k$,
that is $i_s \leq k_s$ for all $s \leq t-1$.  Hence $j \leq k$ and by
assumption also $k \leq j$.  So $\tbinom{i_t p^t + k}{j} \neq 0$
precisely if $k=j$.  In that case, $\tbinom{i_t p^t + k}{j} =
\tbinom{i_t p^t + j}{j} = \tbinom{j}{j} = 1$, $\tbinom{i_t p^t + k}{k}
= \tbinom{k}{k} = 1$, and $\tbinom{m+i_tp^t+j}{j-k} = 1$.
Hence $b({j}) \cdot b({i_t p^t})=b({i_t p^t +j})$.
\end{proof}

\begin{lem} \label{lem13}
We define the degree of the basis element $b(i)$ to be $i$.
Let $1 < n \leq p-1$, then in $S_K(\lambda)$ we have
\[
(b({p^t}))^n= (n!)^2\, b({n \cdot p^t}) + 
\mbox{terms of lower degree.}
\]
\end{lem}

\begin{proof}
This follows by induction on $n$, using the multiplication formula
given in Proposition~\ref{prop2}. More precisely, let $2 \leq c \leq
p-1$, then
\begin{equation*}
b({p^t}) \cdot b({(c-1) p^t}) = \sum_{k=0}^{p^t} 
\tbinom{(c-1)p^t +k}{p^t} \tbinom{(c-1)p^t +k}{k} \tbinom{m+c p^t}{p^t-k}
\; b({(c-1)p^t +k}),
\end{equation*}
where for $k=p^t$ we obtain $\tbinom{(c-1)p^t +k}{p^t} =
\tbinom{cp^t}{p^t} = \tbinom{c}{1}$ and $\tbinom{(c-1)p^t +k}{k} =
\tbinom{cp^t}{p^t}= \tbinom{c}{1}$ and $\tbinom{m+c p^t}{p^t-k} =
\tbinom{m+cp^t}{0}=1$.
Thus the above formula takes the form
\begin{equation}\label{eqn:300}
b({p^t}) \cdot b({(c-1) p^t}) = c^2\, b(c\cdot p^t) 
        + \mbox{\em terms of lower degree.}
\end{equation}
Taking $c=2$ in this formula gives
\begin{equation*}
b({p^t})^2  = 2^2\, b(2\cdot p^t) 
        + \mbox{\em terms of lower degree}
\end{equation*}
and multiplying this through by $b(p^t)$ 
and using Equation \eqref{eqn:300} again yields 
\begin{equation*}
b({p^t})^3  = 3^2 \cdot 2^2\, b(3\cdot p^t) 
        + \mbox{\em terms of lower degree}
\end{equation*}
and so forth.
\end{proof}

\begin{cor}
Let $\lambda=(\lambda_1, \lambda_2)$ be a partition of $r$ and assume
$t$ is such that $p^t \leq \lambda_2 < p^{t+1}$.  Then the algebra
$S_K(\lambda)$ is generated by the elements $b(p^0), b({p^1}), \ldots,
b({p^t})$.
\end{cor}

\begin{proof}
We know already from Lemma~\ref{lem12} a factorisation of
a basis element $b(i)$. Write $i=[i_0, i_1, \ldots]$ $p$-adically. Then
\begin{equation*} \textstyle
b(i) = \prod_{t \geq 0} b({i_t \cdot p^{t}}). 
\end{equation*}
Hence we need to show that the elements $b({ c \cdot p^{t}})$ for $1
\leq c \leq p-1$ are generated by the elements $b({p^t})$. This follows
by induction on $t$ using Lemma~\ref{lem13}. 
\end{proof}

\begin{rk*}
  For odd primes it seems hard to find explicit expressions for
  $b(p^t)^n$ in terms of generators, and arithmetic conditions seem to
  be quite complicated.  For $p=2$ this is done in the next section.
\end{rk*}

\begin{ex} \label{ex:chap3end}
Let $m=0$ and $p=2$. Then $r$ is even and $\lambda=(r/2,r/2)$; in this
case the algebra $S_K(\lambda)$ has dimension $r/2+1$. It is generated
by $b(0), \ldots, b(2^k)$ where $2^k \leq r/2+1 < 2^{k+1}$ subject to
the relations
$$\begin{array}{ll} 
b(2^i)^2=0, & 0\leq i \leq k;\\ 
\prod_{i \in I} b({2^i}) =0, & \mbox{whenever } 
I\subseteq \{ 0, 1, \ldots, k\} \mbox{ and }  
\sum_{i \in I} 2^i \geq r/2+1.
\end{array}
$$ 
It follows that there are no non-zero primitive idempotents except
$1$, and hence $S_K(\lambda)$ is indecomposable; that is, the algebra
is a block.
\end{ex}


\section{The elements $b(i)^2$ are idempotents} \label{sec4} 
From now we assume that the characteristic of the underlying field $K$
is $p=2$.  Then Lemma~\ref{lem12} shows that the basis element $b(i)$
in $S_K(\lambda)$ is equal to the product of the $b(2^t)$ for which
$i_t=1$.  So to understand the multiplication completely we need to
understand the squares of the basis elements $b(2^t)$.

\begin{ex} \label{ex:4.1}
Let $m$ be fixed with $2$-adic expansion $m = [m_0, \ldots, m_t,
\ldots ]$. Suppose $t=0,1$ then we see directly from multiplication
formula (\ref{eqn1a}) that
\begin{equation*}
\begin{aligned}
b(2^0)^2 &= m_0\cdot b(2^0),  \\
b(2^1)^2 &= b(2^1)[ m_1\cdot 1  +  m_0 \cdot b(2^0)].
\end{aligned}
\end{equation*}
So we can write $b(2^1)^2 = b(2^1)(m_1+b(2^0)^2)$. This has the
following generalization.
\end{ex}

\begin{lem}  \label{lem2} 
Suppose $m = [m_0, \ldots, m_t, \ldots ]$ in $2$-adic expansion. Let
$0 \leq v\leq t$ be maximal such that $m_{v-1}=0$. Then
\[
b(2^t)^2 = b(2^t)[ m_t\cdot 1 + \sum_{i=v-1}^{t-1} b(2^i)^2],
\]
setting $b(2^i)=0$ and $m_i=0$ if $i<0$.
\end{lem}

\begin{proof} We make the convention that $m_i=0$ when $i<0$. 
We rewrite the product $b(2^t)^2$ using the multiplication formula
given in Equation (\ref{eqn1a}). 
Note that $\tbinom{2^t+k}{2^t} = \tbinom{2^t+k}{k}$ is zero modulo two when
$k=2^t$, and if $k < 2^t$ it is one modulo two. Moreover for $k<2^t$ we have
\begin{eqnarray} \label{eqn2}
\tbinom{m+2^{t+1}}{2^t-k} \equiv \tbinom{m}{2^t-k} \mod{2}.
\end{eqnarray}
We will change variables using the relation $2^t+k=2^{t+1}-(2^t-k)=
2^{t+1}-l$. Hence -- by Equation (\ref{eqn2}) -- we can rewrite
Equation \eqref{eqn1a} in the form
\begin{eqnarray} \label{eqn*}
b(2^t)^2 &=& \sum_{k=0}^{2^t-1} \tbinom{m}{2^t-k} b(2^t+k) 
= \sum_{l=1}^{2^t}\tbinom{m}{l} b(2^{t+1}-l) 
= b(2^t)[\sum_{l=1}^{2^t}\tbinom{m}{l} b(2^t-l)]. 
\end{eqnarray}
For the last equality note that $2^{t+1}-l = 2^t + (2^t-l)$, and
so for $0\leq 2^t-l < 2^t$ we can factor $b(2^{t+1}-l) =
b(2^t)b(2^t-l)$ by Lemma~\ref{lem12}.
The  term with $l=2^t$ is equal to $m_tb(0) = m_t\cdot 1$.
So we can write
\begin{eqnarray} \label{eqn3}
b(2^t)^2 = b(2^t)[m_t\cdot 1 + \Gamma(t)]
\mbox{\hspace{.5cm} where \hspace{.5cm}} 
\Gamma(t):= \sum_{l=1}^{2^t-1} \tbinom{m}{l}b(2^t-l).
\end{eqnarray}
We will now prove a recursion formula for $\Gamma(t)$. We claim that
\begin{equation}\label{eqn5}
\begin{aligned}
\Gamma(1) &= b(2^0)^2, \\ 
\Gamma(t) &=  b(2^{t-1})^2 +
m_{t-1}\Gamma(t-1) \quad \mbox{ for $t\geq 2$.}
\end{aligned}
\end{equation}

First, $\Gamma(1) = \binom{m}{1}b(2^0) = m_0b(2^0) = b(2^0)^2$, where
the last equality is by Example \ref{ex:4.1}. Suppose that $t\geq 2$,
and split $\Gamma(t)$ into two sums as follows:
\begin{eqnarray*}
\Gamma(t) &=& \sum_{l=1}^{2^{t-1}} {m\choose l}b(2^t-l) +
\sum_{l=2^{t-1}+1}^{2^t-1}{m\choose l}b(2^t-l) 
\mbox{\hspace{.5cm} by definition of $\Gamma(t)$, }\\
&=&  b(2^{t-1})^2 +   \sum_{l=2^{t-1}+1}^{2^t-1}{m\choose l}b(2^t-l) 
\mbox{\hspace{2.1cm} by Equation (\ref{eqn*}), Lemma~\ref{lem12},} \\
&=&  b(2^{t-1})^2 +    m_{t-1} \sum_{r=1}^{2^{t-1}-1} {m\choose r}b(2^{t-1}-r)      \mbox{\hspace{1cm} (the argument is given below),} \\
&=& b(2^{t-1})^2 +  m_{t-1}\Gamma(t-1) 
\mbox{\hspace{3.6cm} by definition of $\Gamma(t-1)$. }
\end{eqnarray*}
For the third equality sign in the latter equation, set $l = 2^{t-1}+
r$ where $1 \leq r \leq 2^{t-1}-1$, and note that
\[
{m\choose 2^{t-1}+r}\equiv m_{t-1}{m\choose r} \mod{2},
\]
and $2^t-l = 2^{t-1}-r$.  Hence the recursion formula for $\Gamma(t)$
claimed in Equation (\ref{eqn5}) is shown.  It in fact implies the
following simpler formula for $\Gamma(t)$:
\begin{equation*}
\Gamma(t) = \sum_{i=v-1}^{t-1} b(2^i)^2
\end{equation*}
where $v$ is as in the statement. Substituting this into (\ref{eqn3})
completes the proof.
\end{proof}

\begin{prop}  \label{cor2}
Let $p=2$.  For $i \geq 0$, the elements $b(2^i)^2$ are idempotent in
$S_K(\lambda)$.  Moreover, if $m_j=0$ for all $j\leq i$ then
$b(2^i)^2=0$.
 \end{prop}

\begin{proof}
This follows from Lemma~\ref{lem2} by induction.
\end{proof}

\section{Analysis of the binomial coefficient $B(m,g)$}

Still keeping $p=2$ fixed, we assume throughout this section that $m$
and $g$ are integers such that the binomial coefficient $B(m,g)$ is
non-zero modulo two.  We need to relate the binomial expansion of $m$ with
that of $B(m,g)$. Note that we have the following depiction of the
binary addition:
\[
\begin{array}{r|cccccc}
m & m_0 &m_1& m_2 & \ldots & m_i & \ldots \\
+2g & 0 & g_0 & g_1 & \ldots & g_{i-1} & \ldots \\
\hline
m+2g & (m+2g)_0 & (m+2g)_1 & (m+2g)_2 & \ldots & (m+2g)_i & \ldots
\end{array}
\]
In this addition, we need to keep track over the `carry overs'.  So
define integers $x_i\geq 0$ such that
\begin{eqnarray} \label{eqn6}
m_i+g_{i-1} + x_{i-1} = (m+2g)_i + 2x_i.
\end{eqnarray}
Thus $x_i$ is the carry over from column $i$ to column $i+1$ in the
binary addition of $m$ and $2g$.  Most important for the proofs later
will be that $(m+2g)_i=1$ implies that $x_i=0$; more precisely we have
the following:

\begin{prop} \label{prop1}
Let $m=[m_0, m_1, \ldots ]$ and $g = [g_0, g_1, \ldots]$ be in binary
expansion.  Assume that $B(m,g)$ is non-zero.  Then $(m+2g)_i + 2x_i <
3$ for all $i$. In particular, if $(m+2g)_i=1$ then $x_i=0$.
\end{prop}

\begin{proof} Certainly $(m+2g)_i + 2x_i \leq 3$. Assume for a
contradiction that this number is equal to three for some $i$. Then
$x_{i-1}=g_{i-1}=1$.  Since $g_{i-1}=1$ we must have that
$(m+2g)_{i-1}=1$ as well, since otherwise the binomial coefficient
$B(m,g)$ would be zero.  But then it follows that $m_{i-1} + g_{i-2} +
x_{i-2} = 3$, and then repeating the argument gives
$m_1+g_0+x_0=3$. This implies $x_0=1$. On the other hand,
$(m+2g)_0=m_0$ and hence $x_0=0$, a contradiction.
\end{proof}

We will later prove some properties by induction. The elements
$e_{m,g}$ are defined as products, and it will be convenient to use
factors of these which are already known to be idempotents.  The basis
for the induction will be the following:

\begin{lem}[Splitting Lemma] \label{lem9}
Let $u$ be a natural number and define 
\[
n:= [m_0, m_1, \ldots, m_u] \mbox{ \hspace{.3cm} and \hspace{.3cm}}
d:= [g_0, g_1, \ldots, g_{u-1}].
\]
Suppose $(m+2g)_u=1$.  Then the binary expansion of $B(n,d)$ equals
the binary expansion of $B(m,g)_{<u}$ extended by one column
$\tbinom{1}{0}$.
In particular if $g_u =0$ then $B(n,d) = B(m,g)_{\leq u}$. 
\end{lem} 
    
\begin{proof} By Proposition~\ref{prop1} we know that $x_u=0$, and by
Equation (\ref{eqn6}) we hence have $m_u +g_{u-1} + x_{u-1} = 1$; the
claim follows.
\end{proof}

\begin{rk*}
The Splitting Lemma shows that when $g_u=0$ then the
element $e_{n,d}$ is a factor of $e_{m,g}$, when written as in the
definition; see Equation \eqref{eqn100}.
\end{rk*}

We will have to use the formula from Lemma~\ref{lem2}. So we need to
know the digits of $B(m,g)$, given the binary expansion of $m$ and of
$g$. We now describe these explicitly.

\begin{lem} \label{lem3} 
Given natural numbers $t$ and $a$.  Suppose $B(m,g)_{\leq \
t+a}$ in binary decomposition is of the form
\begin{eqnarray} \label{eqn7}
B(m,g)_{\leq t+a}= 
\left(\begin{smallmatrix}\ldots & 1&0&\ldots & 0 \cr \ldots &g_t&0&\ldots
& 0\end{smallmatrix}\right). 
\end{eqnarray}
Then we have:
\begin{itemize}
\item[(a)] Suppose $g_t=0$, then $m_{t+1} =  \ldots = m_{t+a}=0$ and
$x_{t+1} = \ldots = x_{t+a}=0$.
\item[(b)] Suppose $g_t=1$, then $m_{t+1} = \ldots = m_{t+a}=1$ and
$x_{t+1} = \ldots = x_{t+a}=1$.
\end{itemize}
\end{lem}

\begin{proof} By Proposition~\ref{prop1} we know that $x_t=0$.  By
Equation (\ref{eqn6}) we have:
\begin{eqnarray*}
m_{t+1} + g_t + 0 &=& 0+2x_{t+1},\\
m_{t+2} + 0 + x_{t+1} &=& 0+2x_{t+2}, \\
\ldots & & \ldots \\
m_{t+a}  + 0  + x_{t+a-1} &=& 0+2x_{t+a}.
\end{eqnarray*}

For (a), assume that $g_t=0$. Then $x_{t+1}=0$ and hence
$m_{t+1}=0$. Now the second equation shows that $x_{t+2}=0$ and hence
$m_{t+2}=0$, and so on.  Part (b) is similar.
\end{proof}

We will have to consider sequences of digits such that $m_i=1$ for
$v\leq i\leq s$ and $m_{v-1}=0$. For these values of $i$ we need to
know the $i$-th columns of $B(m,g)$.

\begin{lem} \label{lem4}  
Suppose column $s$ of $B(m,g)$ is zero but column $s-1$ is non-zero.
Let $u\geq 0$ be minimal such that $(m+2g)_i=1$ for $u\leq i < s$, and
let $0 \leq v \leq s$ be maximal with $m_{v-1}=0$.  Then $v\geq
u$. Moreover:
\begin{itemize}
\item[(a)] If $m_s=0$ then $g_i=0$ for $v-1\leq i\leq s-1$.
\item[(b)] If $m_s=1$ then $g_{s-1}=1$ and $g_i=0$ for $v-1\leq i < s-1$.
\end{itemize}
\end{lem}

\begin{proof} (i) Suppose $m_u=0$ or $u=0$. Then by definition of $v$
we have that $v \geq u$. So assume that $m_u=1$ and $u>0$.  By
definition of $u$ we have that $(m+2g)_u=1$ and $(m+2g)_{u-1}=0$.
Then Equation (\ref{eqn6}) for columns $u$ and $u-1$ together with the
assumptions and Proposition~\ref{prop1} read:
\begin{equation*}
\begin{aligned}
1 + g_{u-1} + x_{u-1} &= 1, \\ 
m_{u-1}+ g_{u-2} + x_{u-2} &= 0+2x_{u-1}.
\end{aligned}
\end{equation*}
So $x_{u-1}=0=g_{u-1}$ which implies that $m_{u-1} + g_{u-2} + x_{u-2}
= 0$ and hence $m_{u-1}=0$. This shows that $u \leq v$.

(ii) For (a) and (b), use Proposition~\ref{prop1} and Equation
\eqref{eqn6} for columns between $v$ and $s-1$. By assumption and (i)
we have that $(m+2g)_i=1=m_i$  for $v \leq i \leq s-1$. This
implies $g_{i-1}=0=x_{i-1}$ for $v \leq i \leq s-1$ and $x_{s-1}=0$.
Then Equation \eqref{eqn6} for column $s$ becomes
\begin{equation*}
m_s+g_{s-1} + 0 = 0+2x_s.  
\end{equation*}
If $m_s=0$ then $x_s=0$ and
$g_{s-1}=0$. On the other hand if $m_s=1$ then $x_s=1$ and
$g_{s-1}=1$.
\end{proof}

\section{The proofs of the Orthogonality Lemma and 
the Idempotent theorem} \label{sec6}

In this section we return to the analysis of the basis $\{ b(i) \}$ of
$S_K(\lambda)$, still under the assumption char$(K) = 2$.  With the
information obtained in the preceding section, we are now in a
position to complete the proof of both the Orthogonality Lemma and the
Idempotent Theorem, stated in Section~\ref{sub:main}.

\subsection{Proof of the Orthogonality Lemma}\label{orthlem}
Suppose the $s$-th column of $B(m,g)$ is zero.  The aim is to show
that $e_{m,g}^2 \cdot b(2^s)^2 =0$.  Recall from Lemma~\ref{lem2} that
$b(2^s)^2 = b(2^s)\psi$ with
\begin{eqnarray} \label{eqn10}
\psi = \psi_{m,s} = m_s + \sum_{i=v-1}^{s-1}b(2^i)^2
\end{eqnarray}
where $0\leq v \leq s$ is maximal such that  $m_{v-1}=0$. We will 
prove that
\begin{eqnarray} \label{eqn101}
(e_{m,g})^2_{<s} \cdot \psi_{m,s} = 0.
\end{eqnarray}
Certainly this then implies the Orthogonality Lemma in Section~\ref{sec2}.
Note that if $s=0$ then $\psi =m_0 =0$ since $(m+2g)_0=m_0=0$. So assume
$s>0$. If all columns before column $s$ are zero then $m_i=0$ for
$i\leq s$ and then $\psi=0$ by Proposition~\ref{cor2}. So assume now
that $w<s$ is such that $(m+2g)_w=1$ and $(m+2g)_i=0$ for
$w+1\leq i\leq s$. We use induction on the number of zero columns
between $w$ and $s$ to prove Equation (\ref{eqn101}).

Suppose column $s-1$ is non-zero. Let
$u\geq 0$ be minimal such that $(m+2g)_i=1$ for $u\leq i < s$. We
apply Lemma~\ref{lem4}, which shows that $v\geq u$.
Moreover, suppose $m_s=0$, then by part (a) of the Lemma we know that
$(e_{m,g})_{<s}$ has factors $(1-b(2^i))$ for $v-1\leq i\leq s-1$.
This gives that $(e_{m,g})^2_{<s} \cdot \psi = 0$ by
Proposition~\ref{cor2}.

Similarly, if $m_s=1$ then part (b) of the Lemma shows that
$(e_{m,g})_{<s}$ has factors $(1-b(2^i))$ for $v-1\leq i < s-1$ and
also a factor $b(2^{s-1})$.  Then the claim follows again from
Proposition~\ref{cor2}, using that $b(2^{s-1})^2\cdot (m_s +
b(2^{s-1})^2) = 0$. This proves the base case of the induction.

For the inductive step, suppose now that column $s-1$ is zero. The
inductive hypothesis states that
\begin{equation*}
(e_{m,g})^2_{< s-1} \cdot \psi_{m,s-1} = 0.
\end{equation*}
If $g_w=0$ then we have by Lemma~\ref{lem3} that $m_i=0$ for $w+1\leq
i\leq s$. Then $v=s$ and we can write
\begin{equation*}
\psi_{m,s} = b(2^{s-1})^2 = \psi_{m,s-1}\cdot b(2^{s-1}),
\end{equation*}
using Lemma~\ref{lem2}.  By the inductive hypothesis we deduce
$(e_{m,g})^2_{<s}\cdot \psi_{m,s} = 0$.
Now suppose $g_w=1$, then by Lemma~\ref{lem3} we know that $m_i=1$ for
$w+1\leq i\leq s$. We rewrite and again use Lemma~\ref{lem2}:
\begin{equation*}
\psi_{m,s} = \psi_{m,s-1} + b(2^{s-1})^2 = \psi_{m,s-1} +
\psi_{m,s-1}\cdot b(2^{s-1}),
\end{equation*}
and again using the inductive hypothesis we have
$(e_{m,g})^2_{<s}\cdot \psi_{m,s} = 0$.  This completes the proof of
\eqref{eqn101}, and hence also the proof of the Orthogonality Lemma.

\subsection{Proof of the Idempotent Theorem}

This will be done by induction on $t$, the largest column label of a
non-zero column in the binary decomposition of $B(m,g)$, which we call
the degree of $e_{m,g}$.  In fact, we will prove the following:

{\bf Claim:} Elements $e_{m,g}$ and $(e_{m,g})_{<t}$ are idempotents.

Assume that $t=0$ then $B(m,g)=\tbinom{1}{0}$.  In particular $m_0=1$
and so $e_{m,g}=(1-b(2^0))$ is idempotent. Also $(e_{m,g})_{<t}=1$ is
idempotent.
We assume the statement holds for all $e_{n,d}$ of degree $<t$. Let
$e_{m,g}$ be of degree $t$ and write $e:= e_{m,g} =P\cdot (1-b(2^t))$
where $P = (e_{m,g})_{< t}$. We have
\begin{equation*}  \label{eqn102}
e^2 = P^2(1-b(2^t)^2) = P^2(1-\psi b(2^t))
\end{equation*}
where $\psi = \psi_{m,g}$ is defined as in Equation (\ref{eqn10}).  We
will show that $P^2\cdot \psi = P^2$, and secondly that $P^2=P$.  This
then implies that $e=e_{m,g}$ is idempotent.
\medskip

(a) We claim that $P^2\cdot \psi = P^2$, that is $P^2(1-\psi)=0$. To
see this, let $\widetilde{m}:= m + 2^t$, then $\widetilde{m}_t=1 +
m_t$ and $\widetilde{m}_i=m_i$ for $i <t$.  Hence $B(\widetilde{m},g)$
differs from $B(m,g)$ in columns $t$ and $t+1$. Therefore
\begin{equation*}
{(e_{m,g})}_{<t} = (e_{\widetilde{m},g})_{<t} = P.
\end{equation*}
Moreover (using $p=2$) we have $\psi_{\widetilde{m}, t} =
1-\psi_{m,t}$. So we get from the Orthogonality Lemma, see Equation
(\ref{eqn101}):
\begin{equation*}
P^2(1-\psi_{m,t}) = (e_{\widetilde{m},g})^2_{<t} \cdot
\psi_{\widetilde{m}, t} = 0.
\end{equation*}

(b) We claim that $P^2=P$. This is clear if $P=1$. So suppose $P>1$,
then there is some $u<t$ maximal such that $(m+2g)_u=1$. If $g_u=0$
then $P=e_{n,d}$ with $d$ and $n$ as in the Splitting
Lemma~\ref{lem9}. Hence by the inductive hypothesis $P$ is idempotent.
If $g_u=1$, then $P=(e_{m,g})_{<u}\cdot b(2^u)$. Define $n$ and $d$ by 
\begin{eqnarray} \label{eqn66}
e_{n,d}=(e_{m,g})_{<u}\cdot (1-b(2^u)).
\end{eqnarray}
By construction $e_{n,d}$ has degree $u<t$ and hence by the inductive
hypothesis we get that $e_{n,d}$ and $(e_{m,g})_{<u}$ are
idempotents. Since the characteristic of the underlying field is two
and by Equation (\ref{eqn66}), we have that $(e_{m,g})_{<t}= P=
(e_{m,g})_{<u} \cdot b(2^u) = e_{n,d} + (e_{m,g})_{<u}$ is
idempotent.  \hfill $\Box$


\section{The correspondence between idempotents and Young modules.}
\label{sec7}

Fix an integer $g\geq 0$ such that $\tbinom{m+2g}{g} \neq 0$. Then we
have for each $r\geq m+2g$ of the right parity a partition $\lambda$ with
$\lambda_1-\lambda_2=m$, and a partition $\mu=(\mu_1, \mu_2)$ with
$\mu_1-\mu_2=m+2g$.  We also have the primitive idempotent $e_{m,g}$
defined in Equation (\ref{eqn100}); and we know that $Y^{\mu}$ is a
direct summand of $M^{\lambda}$. We will now show that in fact
$e_{m,g}$ is the projection of $M^{\lambda}$ corresponding to
$Y^{\mu}$.

\begin{thm} \label{thm2} Let $\lambda, \mu$ be two-part partitions
  of $r$ such that $Y^{\mu}$ is a direct summand of $M^{\lambda}$. Let
  $\lambda_1-\lambda_2=m$, \ $\mu_1-\mu_2 = m+2g$ and
  $g=\lambda_2-\mu_2$. Then the idempotent $e_{m,g}$ of $S_K(\lambda)$
  is the projection onto $Y^{\mu}$.
\end{thm}

The proof of this will take the rest of the section. We use induction
on $r$, starting with the case $\mu_2=0$, that is $\mu= (r, 0)$. Then
the inductive step will be to show that if the theorem is true for
degree $r$ then it is true for degree $r+2$.

To begin the induction we make two observations: 

(i) Suppose that $\mu_2=0$. In the special case when $\lambda=\mu$ we
have $g=0$ and $m=r$. So $\lambda_2=0$ and the algebra $S_K(\lambda)$
has dimension one.  Furthermore, $e_{m,0}=1$ and $M^{\lambda} =
Y^{\lambda}$, so the theorem is trivially true.

(ii) Suppose next that $\mu_2=0$ and $\mu > \lambda$. We have then
$r=\mu_1$ and $\mu_2=0$.  By case (i), we know that $e_{r, 0} \in
S_K(\mu)$ is the projection corresponding to the summand $Y^{\mu}$ of
$M^{\mu}$.  Both idempotents $e_{m,g}$ and $e_{r, 0}$ lie in $S_K(2,
r)$. To show that the summand of $M^{\lambda}$ corresponding to the
projection $e_{m,g}$ is isomorphic to $Y^{\mu}$ we must show that the
idempotents $e_{m,g}$ and $e_{r,0}$ are associated in $S_K(2,r)$.

\begin{prop} 
Under the assumptions in (ii), the idempotents $e_{m,g}$ and $e_{r,0}$
are associated in $S_K(2,r)$. Hence $e_{m,g}M^{\lambda}$ is isomorphic
to $Y^{\mu}$.
\end{prop}

\begin{proof} 
(a) We first simplify the expressions for the two idempotents. Note
that by definition (see Equation \eqref{eqn100}) we have
\begin{eqnarray*}
e_{m,g} 
&=& \prod_{u \in J_{m,g}} b(2^u)\cdot \prod_{u \in I_{m,g}} (1 -b(2^u)) 
\mbox{\hspace{1cm} by Equation (\ref{eqn100}),} \\
&=& b(g) \cdot \prod_{u \in I_{m,g}} (1 -b(2^u)) 
\mbox{\hspace{2.3cm} by Lemma~\ref{lem12}, }\\
&=& b(g) \cdot (1 \pm \mbox{sum of products of $b(i)$'s}\,\, ) \\
&=& b(g)
\end{eqnarray*}
To get the last equality note that $b(g) \cdot b(i)=0$ as the algebra
$S_K(\lambda)$ has basis $\{b(0), b(1), \ldots, b(g)\}$ and by using
Lemma~\ref{lem12}.  Moreover, as $M^{(r,0)}=Y^{(r,0)}$, we have
$e_{r,0}=1_{(r,0)}$.

(b) Let $\alpha=(1,-1)$ and recall from~\cite{DG:PSA}, Theorem 2.4,
that for any composition  $\nu$ we have
\begin{eqnarray*} 
e \cdot 1_\nu = \left\{
\begin{array}{ll}
1_{\nu+\alpha} \cdot e & \mbox{ if $\nu +\alpha$ is a composition,} \\
0 & \mbox{otherwise, }
\end{array}
\right.
\end{eqnarray*}
and 
\begin{eqnarray*} 
f \cdot 1_\nu = \left\{
\begin{array}{ll}
1_{\nu-\alpha} \cdot e & \mbox{ if $\nu -\alpha$ is a composition,} \\
0 & \mbox{otherwise, }
\end{array}
\right.
\end{eqnarray*}
Moreover, by~\cite{DG:PSA}, Proposition 4.3 we have that $H_i \cdot
1_{\lambda}= \lambda_i \cdot 1_{\lambda}$ for $i=1,2$, and recall that
$h=H_1-H_2$.
These formulas imply that $e \cdot 1_{(r,0)}=0$ as $(r,0)+\alpha$ is
not a composition. Moreover, with $\lambda=(g+m,g)$ a partition of
$r=m+2g$ we have
\begin{equation*}
e^{(g)} \cdot 1_{\lambda} = 1_{(r,0)} \cdot e^{(g)},  \qquad
  1_{(r,0)}\cdot f^{(g)} = f^{(g)} \cdot 1_{\lambda}, \qquad
  \tbinom{h}{g} \cdot 1_{(r,0)}= \tbinom{r}{g} \cdot 1_{(r,0)}.
\end{equation*}

(c) We next give elements $u$ and $v$ in the Schur algebra $S_K(2,r)$
such that $e_{m,g}=uv$ and $e_{r,0}=vu$, proving that the two
idempotents are associated.  More precisely, let
\begin{equation*}
u=1_\lambda f^{(g)} 1_{(r,0)}  \qquad  \text{and} \qquad
v= 1_{(r,0)} e^{(g)} 1_\lambda. 
\end{equation*}
Then by repeated use of the equations in (b) we have
\begin{eqnarray*}
  u \cdot v &=& 1_\lambda \, f^{(g)} \, 1_{(r,0)} \, e^{(g)} \, 1_\lambda 
  =  1_\lambda \, f^{(g)} \, e^{(g)} \, 1_{\lambda} =b(g)
\end{eqnarray*}
and
\begin{eqnarray*}
  v \cdot u &=& 1_{(r,0)} \, e^{(g)} \,1_\lambda \,f^{(g)} \,
1_{(r,0)} \\ &=& 1_{(r,0)} \,e^{(g)} \, f^{(g)} \,1_{(r,0)} \\ 
&=&
1_{(r,0)} \cdot \, [\sum_{j=0}^{g} f^{(g-j)}
\tbinom{h-2g+2j}{j}e^{(g-j)}] \,\cdot 1_{(r,0)} \\ 
&=& 1_{(r,0)} \cdot
[f^{(0)} \tbinom{h}{g} e^{(0)}] \cdot 1_{(r,0)} \\ 
&=& \tbinom{r}{g}
\cdot 1_{(r,0)} =B(m,g) \cdot 1_{(r,0)} = 1_{(r,0)}
\end{eqnarray*}
modulo two.  Hence $e_{m,g}=b(g)$ and $e_{r,0}= 1_{(r,0)}$ are
associated.
\end{proof}

We hence have shown that whenever $\mu=(r,0)$ then the claim made in
Theorem~\ref{thm2} is true.  
Now it remains to deal with the inductive step. We assume that
Theorem~\ref{thm2} holds in degree $r$ and show it holds in degree
$r+2$. Clearly any pair of partitions $\tilde{\lambda} < \tilde{\mu}$
in degree $r+2$ with $\tilde{\mu}_2>0$ which satisfies the assumptions
of Theorem~\ref{thm2}, is obtained from partitions $\lambda > \mu$ in
degree $r$ as $\tilde{\mu} = \mu+(1,1)$ and
$\tilde{\lambda}=\lambda+(1,1)$. We hence do the induction step by
comparing $M^{\lambda}$ and $M^{\lambda + (1^2)}$.  To do so, we will
first analyze more closely how the hyperalgebra actions on $E^{\otimes
r}$ and $E^{\otimes r+2}$ are related.
We fix a basis $\{v_1, v_2\}$ of the $K$-vector space $E$. We write briefly
$v_{\di}$ for the tensor product $v_{i_1} \otimes v_{i_2}
\otimes \ldots \otimes v_{i_r}$, with  $\di$ the  multi-index 
$\di=(i_1, \ldots,
i_r)$.  Define the linear map
\begin{equation*}
j: E^{\otimes r} \longrightarrow E^{\otimes r+2} \mbox{\,\, by \,\,} x
\mapsto (v_1 \otimes v_2 - v_2 \otimes v_1) \otimes x.
\end{equation*}

Recall that both tensor powers  are modules for the hyperalgebra $\U_K =
\U(\gl_2)_Z\otimes K$. The map $j$ commutes with the action of the
divided powers $e^{(a)}, \ f^{(a)} \in \U_K$: this is easy to
see, noting that
the map $j$ is tensoring with
$\bigwedge^2E$, which is trivial under the action of $e$ and $f$.
Now we restrict $j$ to $M^{\lambda}$; it takes
$M^{\lambda}$ to $M^{\lambda + (1^2)}$. Since the products
$f^{(a)}e^{(a)}$ lie in the zero weight space of $\U_K$, they preserve
$M^{\lambda}$ and $M^{\lambda + (1^2)}$. The idempotents $1_{\lambda}$
and $1_{\lambda + (1^2)}$ are the projections onto these spaces, and
it follows that $j$ intertwines the actions of elements $b(a)$ on
$M^{\lambda}$ and on $M^{\lambda + (1^2)}$.  In particular this implies
\begin{eqnarray} \label{eqn14}
j(e_{m,g}x) = e_{m,g}j(x), \qquad \mbox{ for all $x\in M^{\lambda}$.}
\end{eqnarray}
The following proposition completes the proof of Theorem~\ref{thm2}.

\begin{prop} 
Suppose $e_{m,g}$ is the projection on $M^{\lambda}$ corresponding
to $Y^{\mu}$. Then $e_{m,g}$  on $M^{\lambda + (1^2)}$ is the 
projection corresponding to $Y^{\mu + (1^2)}$. 
\end{prop}

\begin{proof} 
We may assume $m\neq 0$; the case $m=0$ is understood, see
Example~\ref{ex:chap3end}.  We know that the Specht module $S^{\mu}$
is a submodule of $Y^{\mu}$. Furthermore, ${\rm Hom
}_{K\Sigma_r}(S^{\mu}, M^{\lambda})$ is one-dimensional (see
\cite[13.13]{james}). So $M^{\lambda}$ has a unique submodule
isomorphic to $S^{\mu}$, which is contained in $Y^{\mu}$. Similarly
$M^{\lambda + (1^2)}$ has a unique submodule isomorphic to
$S^{\mu+(1^2)}$ and it is contained in $Y^{\mu+(1^2)}$. Since the
elements $e_{m,g}$ are projections onto a Young module, it suffices to
show the following:
\begin{equation*}
\text{
If $e_{m,g}(S^{\mu})\neq 0$ in $M^{\lambda}$ then $e_{m,g}(S^{\mu +
(1^2)})\neq0$ in $M^{\lambda + (1^2)}$. 
}
\end{equation*}
To do so we use polytabloids, that is the standard generators for
Specht modules, see James~\cite{james}, Chapter 4.  We start with
standard tableaux of shapes $\mu$ and $\mu+(1^2)$ respectively, the
two rows of which are filled in as follows:
\begin{eqnarray*}
   \setlength{\unitlength}{1.1mm}
   \begin{picture}(110,8)(0,0)
   \put(0,3){$t_1=$}
   \put(7,4.5){$3\,5\,...\,(2u-1) \,(2u+1)\,...\,(r+2)$}
   \put(7,1.5){$4\,6\,...\,(2u) \,$}
   \put(63,3){$t_2=$}
   \put(70,4.5){$1\,3\,...\,(2u-1) \,(2u+1)\,...\,(r+2)$}
   \put(70,1.5){$2\,4\,...\,(2u) \,.$}
   \end{picture}
   \end{eqnarray*}
Here $u=\mu_2+1$. Let $R_{t_i}$ be the row stabilizer of $t_i$, and
$C_{t_i}$ the column stabilizer of $t_i$.  To write down the
polytabloid generating $S^{\mu}$ in this setup, we must start with an
appropriate element $\omega_1 \in M^{\lambda}$ which is fixed by all
elements of $R_{t_1}$. Then the corresponding 'polytabloid' is
\begin{equation*}
\varepsilon_{t_1} =  \omega_1\{C_{t_1}\}^-     
\end{equation*}
where $\{ C_{t_1}\}^-$ is the alternating sum over all elements in 
$C_{t_1}$. We can take
\begin{equation*}
\omega_1 = \sum v_{\di}  
\end{equation*}
summing over all $\di$ such that $i_\rho=2$ for $\rho$ in the second
row of $t_1$, and all other $i_\rho \in \{1,2\}$ such that the weight
of $\di$ is $\lambda$. Note that $\lambda_2\geq \mu_2$, so such a
$\di$ exists.  (When $\lambda=\mu$ then $\omega_1$ consists of just one
basis vector.)  Similarly one defines the Specht module generator
$\varepsilon_{t_2}$ from $t_2$.
Explicitly,
\begin{equation*}
\{ C_{t_1}\}^- = (1- (3,4))(1-(5,6))\ldots (1-(2u-1,2u))
\end{equation*}
This shows that $\omega_1\{C_{t_1}\}^- =
\tilde{\omega}_1\{C_{t_1}\}^-$ where $\tilde{\omega}_1$ is the sum
over all $v_{\di}$ such that $i_{2\rho+1}=1$ and $i_{2\rho+2}=2$ for
$1\leq \rho < u$ (and $i_\rho \in \{1,2\}$ otherwise such that the
weight of $\di$ is $\lambda$).  We next apply the map $j$ to
$\varepsilon_{t_1}$:
$$
j(\varepsilon_{t_1}) =  
(v_1\otimes v_2 - v_2\otimes v_1)\otimes \varepsilon_{t_1}
\ 
= \  (v_1\otimes v_2 \otimes \tilde{\omega}_1(1-(1,2))\cdot \{C_{t_1}\}^-.
$$
Now, $(1-(1,2)\{C_{t_1}\}^- = \{ C_{t_2}\}^-$ and $v_1\otimes v_2
\otimes \tilde{\omega}_1 = \tilde{\omega}_2$.  This shows that $j$
takes $\varepsilon_{t_1}$ precisely to $\varepsilon_{t_2}$.

We can now complete the inductive step of the proof. Suppose
$e_{m,g}(S^{\mu})\neq 0$, then $e_{m,g}(\varepsilon_{t_1})\neq 0$
since $\varepsilon_{t_1}$ is a generator of the Specht module (and
$e_{m,g}$ is a homomorphism). Then also $j\circ
e_{m,g}(\varepsilon_{t_1})\neq 0$ since $j$ is one-to-one. Hence by
Equation (\ref{eqn14}),
\begin{equation*}
0 \neq j (e_{m,g}(\varepsilon_{t_1})) = e_{m,g}\circ j(\varepsilon_{t_1})
= e_{m,g}(\varepsilon_{t_2}).
\end{equation*}
Hence $e_{m,g}(S^{\mu+(1^2)})\neq 0$, as required. 
\end{proof}


\nocite{*} \bibliographystyle{amsplain}

\providecommand{\bysame}{\leavevmode\hbox to3em{\hrulefill}\thinspace}
\providecommand{\MR}{\relax\ifhmode\unskip\space\fi MR }
\providecommand{\MRhref}[2]{%
  \href{http://www.ams.org/mathscinet-getitem?mr=#1}{#2}
}
\providecommand{\href}[2]{#2}


\bigskip
\footnotesize
\begin{tabular}{lcr}
  Stephen Doty & Karin Erdmann & Anne Henke \\
  Mathematics and Statistics \phantom{xxxxxx}& 
  \phantom{xxx} Mathematical Institute \phantom{xxx} & 
  \phantom{xxxxxx} Mathematical Institute \\
  Loyola University Chicago & University of Oxford & University of
  Oxford \\
  Chicago, Illinois 60626 USA & OX1 3LB Oxford, UK & OX1 3LB Oxford, UK\\
  {\tt doty@math.luc.edu} & {\tt erdmann@maths.ox.ac.uk} & {\tt
    henke@maths.ox.ac.uk} \\
\end{tabular}

\end{document}